**Using the Blackwell-Cover Guessing Strategy to Separate Unobserved Bernoulli Trials**

By James D. Stein

Department of Mathematics and Statistics, California State University (Long Beach)

**Abstract –** Blackwell and Cover both proposed the same method for guessing the larger of a pair of two numbers when one can only observe one – the other of the two numbers remains unobserved. As originally stated, the procedure ended with a computation showing that the probability of correctly guessing the larger number was greater than ½.

In this paper, we investigate what we can say when about the number in the pair that was not observed we employ the Blackwell-Cover procedure. Let $\alpha_1$ and $\alpha_2$ be two different numbers in (0,1]. Let G be the uniform distribution on $[0, \alpha_1]$ and H be the uniform distribution on $[0, \alpha_2]$. We think of G and H as being represented by the heads probabilities of coins.

We use a variant of the Blackwell-Cover strategy ro take pairs of coins chosen from distributions created from G and H such that one of the coins is flipped, and as a result of that flip the unobserved (unflipped) coin is placed in one of two sets A and B. We show that there are situations in which this can be done so that the mean of the heads probabilities of the unobserved coins in A is less than the mean of the heads probabilities of the unobserved coins in B. We show that there are situations in which recursive application of this process results in a set consisting almost completely of coins from one of the two distributions G and H, and another set consisting almost completely of coins from the other distribution.

## Introduction

A well-known technique involving using a random variable, ascribed variously to both Blackwell ([1]) and Cover ([2]) allows one to guess which of two numbers is larger with success probability greater than ½, Previous papers ([3], [4]) have investigated what happens when the numbers come from different probability distributions. In such situations, one needs to decide how to evaluate the success of the prediction when the two numbers are equal.

The two obvious methods of doing this are either to not count the guess at all, or to count it as half-correct and half-incorrect. Either method will not change whether the overall strategy being invoked is successful more than half the time, but it will affect the computation of the success probability of the overall strategy.

The Blackwell-Cover strategy can be used to examine items from a set two at a time, placing the unobserved item of the pair in set A if it is guessed to be smaller and in set B if it is guessed to be larger. In this case, what happens when the two items have the same value can affect whether the overall strategy is successful more than half the time – because the Blackwell-Cover strategy does the wrong thing when applied to a pair of equal items.

Suppose that the two items are coins with heads probabilities L and S, where L > S. The Blackwell-Cover strategy says to choose one of the two coins equiprobably and flip it. If it lands heads place the unobserved coin in set A, and if it lands tails, place the unobserved coin in set B. The idea here is to create two sets A and B such that the mean expected heads probability of the coins in A is lower than the mean of the original collection, and the mean expected heads probability of the coins in B is higher than the mean of the original collection.

In this case, the strategy works well. Half the time we will flip the coin with heads probability L. It will land heads with probability L, placing S in A. Half the time we will flip the coin with heads probability S. It will land heads with probability S, placing L in A. The expected number of heads in A from this flip is ½ LS. The expected number of coins placed in A is ½ (L+S), and the expected heads probability from these coins is LS. So the mean expected heads probability of the coins in A is LS/(1/2(L+S)), and it is easy to show that this is less than (L+S)/2, as the inequality $LS/(1/2(L+S)) < (L+S)/2$ is equivalent through cross-multiplication and simplification to the inequality $4LS < (L+S)^2$, which is easily seen to be equivalent to $(L-S)^2 > 0$. A parallel computation shows that the mean heads probabilities of the coins in B is greater than (L+S)/2.

However, if L=S, this strategy does the wrong thing. If $L > ½$, a coin with heads probability greater than ½ will be placed in A with probability $L > 1/2$. – not what one wants to see happen. Similarly, if $L < ½$, a coin with heads probability less than ½ will be placed in B with probability $1-L > 1/2$ – also not what one wants to see happen.

In this paper, we examine two situations in which we can avoid this counterproductive scenario. The first is analogous to the prediction situation where we avoid counting the guess when the two numbers are equal, and the second is analogous to the prediction situation where we count the guess as half-correct and half-incorrect.

**Section I - The Experimental Setup**

In what follows, the results are only interesting in the case that the individual conducting the experiment (flipping and assigning coins) does not know whether a particular coin comes from G or H, as knowing this would render the experiment trivial and uninteresting. The individual conducting the experiment is simply presented with two identical-looking coins. One coin is labeled “Flip Me” and the other is unlabeled.  The individual conducting the experiment flips the coin labeled “Flip Me”.  If it lands heads, he assigns the other coin to set A, and if it lands tails, he assigns the other coin to set B.

The individual or entity preparing the two coins to be given to the individual conducting the experiment has one bag of coins labeled G, which consists of coins with uniformly distributed heads probabilities on $[0,\alpha_1]$, and one bag of coins labeled H, which consists of coins with uniformly distributed heads probabilities on $[0,\alpha_2]$.  He also has a spinner which lands in the white region with probability $1/(1+\beta)$, where $\beta > 0$.  It lands in the red region with probability $\beta/(1+\beta)$.

The individual or entity preparing the two coins follows one of two scenarios.

**Scenario 1** – The spinner is spun once.  If it lands in the white region, a coin is chosen randomly (uniformly) from G and is labeled “Flip Me”.  The other coin is chosen randomly (uniformly) from H.  If the spinner lands in the red region, the “Flip Me” coin is chosen randomly from H and the other coin from G.

**Scenario 2** – The spinner is spun twice. There are three possibilities.

(1) – Both spins land in the white region.  In this case, a fair coin is labeled “Flip Me” and the other coin is chosen randomly from G.

(2) Both spins land in the red region.  In this case, a fair coin is labeled “Flip Me” and the other coin is chosen randomly from H.

(3) One spin lands in the white region and one in the red region. In this case, one coin is chosen randomly from G and one coin is chosen randomly from H. A fair coin is flipped to determine which of these two coins is labeled "Flip Me".

In Scenario 1, every pair of coins consists of one each from both G and H. The ratio of coins from H that were flipped (examined) to coins from G that were flipped (examined) is $\beta$. The unseen coins therefore appear in the reciprocal ratio.

In Scenario 2, the pairs of coins occur with independent probabilities. The probability that the two coins are both from G is $1/(1+\beta)^2$, the probability that both coins are from H is $\beta^2/(1+\beta)^2$, and the probability that one coin is from G and the other from H is $2\beta/(1+\beta)^2$.

In what follows, we let the letter X denotes the eXamined (flipped) coin and the letter U the Unexamined (unflipped) coin.

**An Intuitive Description of Scenarios 1 and 2**

We think of the process of preparation of the two coins as a distribution D constructed by thoroughly mixing 1 bag of coins with heads probabilities uniformly distributed in $[0,\alpha_1]$ (the distribution G) with $\beta$ bags of coins with heads probabilities uniformly distributed in $[0,\alpha_2]$ (the distribution H). The coins are colored to denote which distribution, G or H, they come from (again, these colors are for the purpose of the reader or the preparer to distinguish between G and H and are not seen by the individual conducting the experiment). In both scenarios, first pick a coin from the bag with all the coins in it. In Scenario 1, label this coin "Flip Me", and then pick the second coin, which is to be given to the conductor of the experiment, that is of the opposite color randomly from the bag. In Scenario 2, pick another coin randomly from the bag. If the coins are of the same color, label a fair coin "Flip Me" and give it and one of the two coins to the conductor of the experiment. If the coins are of opposite color, flip a fair coin to determine which of the two coins is to be labeled "Flip Me" and give both coins to the conductor of the experiment.

**Section II – Computations for Scenario 1**

**Theorem 1 (Results for Scenario 1)**

$P(U \to A) = \frac{(\alpha_1+\beta\alpha_2)}{2(1+\beta)}$ $\qquad$ $m(U \to A) = \frac{\alpha_1\alpha_2(1+\beta)}{2(\alpha_1+\beta\alpha_2)}$

$P(U \to B) = \frac{2(1+\beta)-(\alpha_1+\beta\alpha_2)}{2(1+\beta)}$ $\qquad$ $m(U \to B) = \frac{2(\alpha_2+\beta\alpha_1)-(1+\beta)\alpha_1\alpha_2}{4(1+\beta)-2(\alpha_1+\beta\alpha_2)}$

**Proof:** In this case, one coin of the selected pair will be from distribution G and one will be from distribution H.

The examined coin comes from distribution G with probability $1/(1+\beta)$ (this is how $\beta$ manifests itself, in the probability that the examined coin will come from distribution H as opposed to distribution G in the ratio $\beta$ to 1). With probability q, that coin will flip heads, and $U \to A$. The combined probability of distribution G being selected and the coin flipping heads is $q/(1+\beta)$, and the expected number of unseen heads sent to A from this selection is $pq/(1+\beta)$.

With probability 1-q, that coin will flip tails, and $U \to B$. The combined probability of distribution G being selected and the coin flipping tails is $(1-q)/(1+\beta)$, and the expected number of unseen heads sent to B from this selection is $p(1-q)/(1+\beta)$.

The examined coin comes from distribution H with probability $\beta/(1+\beta)$. With probability p, that coin will flip heads, and $U \to A$. The combined probability of distribution H being selected and the coin flipping heads is $p\beta/(1+\beta)$, and the expected number of unseen heads sent to A from this selection is $pq\beta/(1+\beta)$.

With probability 1-p, that coin will flip tails, and $U \rightarrow B$. The combined probability of distribution H being selected and the coin flipping tails is $(1-p)\beta/(1+\beta)$, and the expected number of unseen heads sent to B from this selection is $(1-p)q\beta/(1+\beta)$.

Recalling that q is a random variable from the uniform distribution on $[0,\alpha_1]$ and p is a random variable from the uniform distribution on $[0,\alpha_2]$ enables the following computations.

$$\mathrm{P(U{\rightarrow}A)} = \frac{1}{\alpha_1\alpha_2(1+\beta)}\int_0^{\alpha_2}\int_0^{\alpha_1}(\beta p + q)\,dq\,dp = \frac{(\alpha_1+\beta\alpha_2)}{2(1+\beta)}$$

Notice that this is also the mean heads probability for distribution D.

The expected number of heads in A, EH(A), is computed by evaluating the following integral.

$$\mathrm{EH(A)} = \frac{1}{\alpha_1\alpha_2}\int_0^{\alpha_2}\int_0^{\alpha_1} pq\,dq\,dp = \frac{\alpha_1\alpha_2}{4}$$

So

$$\mathrm{m(U{\rightarrow}A)} = \frac{EH(A)}{\mathrm{P(U{\rightarrow}A)}} = \frac{\alpha_1\alpha_2(1+\beta)}{2(\alpha_1+\beta\alpha_2)}$$

We also have

$$\mathrm{P(U{\rightarrow}B)} = 1 - \mathrm{P(U \rightarrow A)} = \frac{2(1+\beta)-(\alpha_1+\beta\alpha_2)}{2(1+\beta)}$$

$$\mathrm{EH(B)} = \frac{1}{\alpha_1\alpha_2(1+\beta)}\int_0^{\alpha_2}\int_0^{\alpha_1}(p(1-q)+\beta q(1-p))\,dq\,dp = \frac{(\alpha_2+\beta\alpha_1)}{2(1+\beta)} - \frac{\alpha_1\alpha_2}{4}$$

$$\mathrm{m(U{\rightarrow}B)} = \frac{EH(B)}{\mathrm{P(U{\rightarrow}B)}} = \frac{2(\alpha_2+\beta\alpha_1)-(1+\beta)\alpha_1\alpha_2}{4(1+\beta)-2(\alpha_1+\beta\alpha_2)} \quad \blacksquare$$

It is not surprising that $\lim_{\beta \to 0} \mathrm{m}(\mathrm{U} \to \mathrm{A}) = \frac{\alpha_2}{2}$ and $\lim_{\beta \to \infty} \mathrm{m}(\mathrm{U} \to \mathrm{A}) = \frac{\alpha_1}{2}$, as when $\beta \to 0$ D consist mostly of elements from G. Therefore, in a pair of the form {G,H}, X will be from G most of the time, and U from H. The reverse happens as $\beta \to \infty$.

Notice that if $\alpha_1 = \alpha_2$, $\mathrm{m}(\mathrm{U}\to\mathrm{A}) = \alpha_1/2$, the mean of all the coins in D.

It is also not surprising that $\lim_{\beta \to 0} \mathrm{m}(\mathrm{U} \to \mathrm{B}) = \frac{\alpha_2}{2}$ and $\lim_{\beta \to \infty} \mathrm{m}(\mathrm{U} \to \mathrm{B}) = \frac{\alpha_1}{2}$. And again, if $\alpha_1 = \alpha_2$, $\mathrm{m}(\mathrm{U}\to\mathrm{B}) = \alpha_1/2$

**Cor. 1.1** – In Scenario 1, if $\alpha_1 \neq \alpha_2$, $\mathrm{m}(\mathrm{U}\to\mathrm{A}) < \mathrm{m}(\mathrm{U}\to\mathrm{B})$.

**Proof:** This requires showing that

$$\frac{\alpha_1\alpha_2}{2(\alpha_1 + \beta\alpha_2)} < \frac{2(\alpha_2 + \beta\alpha_1) - (1+\beta)\alpha_1\alpha_2}{4(1+\beta) - 2(\alpha_1 + \beta\alpha_2)}$$

By cross multiplication, this is equivalent to showing

$$4\alpha_1\alpha_2(1 + \beta) - 2(\alpha_1 + \beta\alpha_2)\,\alpha_1\alpha_2 < 4(\alpha_1 + \beta\alpha_2)\,(\alpha_2 + \beta\alpha_1) - 2\,\alpha_1\alpha_2(1 + \beta)\,(\alpha_1 + \beta\alpha_2)$$

The term $-2(\alpha_1 + \beta\alpha_2)\,\alpha_1\alpha_2$, which appears on the left side, also appears on the right side when we distribute the last term by multiplying by $1 + \beta$. When we cancel this term and expand, we are left with showing

$$4\alpha_1\alpha_2 + 4\alpha_1\alpha_2\beta < 4\alpha_1\alpha_2 + 4\alpha_1^2\beta + 4\alpha_2^2\beta + 4\alpha_1\alpha_2\beta^2 - 2\alpha_1^2\alpha_2\beta - 2\alpha_1\alpha_2^2\beta^2$$

Cancel $4\alpha_1\alpha_2$ from both sides and divide both sides by $2\beta$. It remains to show

$$2\alpha_1\alpha_2 < 2\alpha_1^2 + 2\alpha_2^2 + 2\alpha_1\alpha_2\beta - \alpha_1^2\alpha_2 - 2\alpha_1\alpha_2^2\beta$$

This is equivalent to showing

$$0 < 2\alpha_1^2 + 2\alpha_2^2 - \alpha_1\alpha_2(2 + a_1) + \alpha_1\alpha_2\beta(2 - \alpha_2)$$

But $2\alpha_1^2 + 2\alpha_2^2 - \alpha_1\alpha_2(2 + a_1) + \alpha_1\alpha_2\beta(2 - \alpha_2) > 2(\alpha_1 - \alpha_2)^2 + \alpha_1\alpha_2\beta(2 - \alpha_2) > 0$. █

The following simple example shows that although m(U→A) < m(U→B), the overall mean m need not necessarily be between these two values.

**Ex. 1** – We assume $\alpha_1 = 1$, $\alpha_2 = 0.5$

| **β** | **m** | **P(U→A)** | **m(U→A)** | **P(U→B)** | **m(U→B)** |
|---|---|---|---|---|---|
| 0.5 | 0.4167 | 0.4167 | 0.3 | 0.5833 | 0.3571 |
| 1 | 0.375 | 0.375 | 0.3333 | 0.625 | 0.4 |
| 2 | 0.3333 | 0.3333 | 0.375 | 0.6667 | 0.4375 |

We remark in passing that it is simple to obtain exact algebraic expressions for when m(U→A) < m < m(U→B) in terms of the values $\alpha_1$, $\alpha_2$ and β should that prove desirable.

**The Ratio β and the Probability That an Element from G Will Be Examined**

We can also describe the ratio β in terms of the probability $P_G$ that, given a pair of elements with one from each of the distributions G and H, we select the element from G for examination. In this case, $P_G = 1/(\beta+1)$, and the formulas for assignment probabilities and means become

$$P(U\to A) = \frac{P_G\alpha_1+(1-P_G)\alpha_2}{2} \qquad m(U\to A) = \frac{P_G\alpha_1\alpha_2}{2(P_G\alpha_1+(1-P_G)\alpha_2)}$$

$$P(U\to B) = 1- \frac{P_G\alpha_1+(1-P_G)\alpha_2}{2} \qquad m(U\to B) = \frac{2(P_G\alpha_2+(1-P_G)\alpha_1)- \alpha_1\alpha_2}{4-2(P_G\alpha_1+(1-P_G)\alpha_2)}$$

**Imputed Values of the Proportionality Constant β for A and B**

Let $A_G = A \cap G$ and $A_H = A \cap H$. The expected number of heads for a coin randomly chosen from $A_G$ is $\alpha_1/2$, and the expected number of heads for a coin randomly chosen from $A_H$ is $\alpha_2/2$ Since the expected number of heads for a coin randomly chosen from A is m(U→A), we let $\beta_A$ denote the ratio of the number of coins chosen from $A_H$ to the number of coins chosen from $A_G$ in order for a coin chosen according to this procedure to have m(U→A) as the expected number of heads. This results in the equation

$$m(U\rightarrow A), = \frac{1\times\frac{\alpha_1}{2}+\beta_A\times\frac{\alpha_1}{2}}{1+\beta_A}$$

Inserting the known value $m(U\rightarrow A) = \frac{\alpha_1\alpha_2(1+\beta)}{2(\alpha_1+\beta\alpha_2)}$ and solving for $\beta_A$ yields

$$\beta_A = \frac{\alpha_1}{\alpha_2}\frac{1}{\beta}$$

Defining $\beta_B$ similarly for m(U→B) and performing the analogous computation yields

$$\beta_B = \frac{2-\alpha_1}{2-\alpha_2}\frac{1}{\beta}$$

This helps us explore the possibility of regarding A and B as initial distributions D and performing the Blackwell-Cover separation algorithm on each of those. When we perform Blackwell-Cover separation on A, we denote by AA the set of unseen coins with the lower mean expected heads probability and AB the set of unseen coins with the higher mean expected heads probability, BA and BB are defined similarly for the set B.

Let $a = \alpha_1/\alpha_2$ and $b = (2 - \alpha_1)/(2 - \alpha_2)$. Then $\beta_A = a/\beta$ and $\beta_B = b/\beta$. Observe that

$$\beta_{AA} = a/\ \beta_A = a/(a/\beta) = \beta$$

Similarly, $\beta_{BB} = \beta$.

**The Semigroup W(A,B) of Words in the Letters A and B, and Associated Betas**

Suppose that e represents the empty word in W(A,B), and $\beta e = \beta$. We can continue inductively to define distributions of unseen coins by the following scheme. Given a word w $\epsilon$ W(A,b), it corresponds to a distribution in which unseen coins from the uniform distribution on $[0,\alpha_1]$ are mixed with unseen coins from the uniform distribution on $[0,\alpha_2]$ so that a coin drawn from this distribution will be from $[0,\alpha_1]$ with probability $1/(1+\beta_w)$. We can think of this as mixing coins from the uniform distribution on $[0,\alpha_1]$ with coins from the unform distribution on $[0,\alpha_2]$ in the ratio 1 to $\beta_w$.

The Blackwell-Cover separation procedure produces two sets wA and wB such that wA can be regarded as mixing coins from the uniform distribution on $[0,\alpha_1]$ with coins from the unform distribution on $[0,\alpha_2]$ in the ratio 1 to $\beta_{wA}$, and similarly for wB. Using the variables a and b as defined in the previous section, we have

$$\beta_{wA} = a/\beta_w \qquad\qquad \beta_{wB} = b/\ \beta_w$$

along with the relations

$$\beta_{wAA} = \beta_{wBB} = \beta_w$$

However, if we alternate the letters A and B, the situation is more complicated and leads to an intriguing result.

**Cor. 1.2 -** $\beta_{(AB)^n} = \boldsymbol{\beta}(\frac{b}{a})^n \quad \beta_{(AB)^nA} = \frac{a}{\beta}(\frac{a}{b})^n \quad$ **for n = 1, 2, ...**

**Proof:** This is easily established by induction. For n=1, $\beta_{AB} = b/\beta_A = b/(a/\beta) = \beta b/a$. We also have $\beta_{ABA} = a/\beta_{AB} = a/(\beta b/a) = (a/\beta)(a/b)$.

Assume the above formulas for n. We then have $\beta_{(AB)^{n+1}} = \frac{b}{\beta_{(AB)^n A}} = \boldsymbol{\beta}(\frac{\boldsymbol{b}}{\boldsymbol{a}})^{\boldsymbol{n+1}}$ and also

$\beta_{(AB)^{n+1}A} = \frac{a}{\beta_{(AB)^{n+1}}} = \frac{\boldsymbol{a}}{\boldsymbol{\beta}}(\frac{\boldsymbol{a}}{\boldsymbol{b}})^{\boldsymbol{n+1}}$. ■

Notice that if $\alpha_1 \neq \alpha_2$, then either a > b or b > a. In either case, one of the fractions a/b and b/a is greater than 1, and the other is less than 1. So the imputed betas for the two sequences $(AB)^n$ and $(AB)^n A$ converge to 0 and $+\infty$. By the remarks following the proof of Theorem 1, the means of the associated sets converge to the means of G and H. We can interpret this as being able to retrieve subsets of unobserved coins in which almost all the coins come from G, or almost all the coins come from H – without ever having flipped any of those coins.

## Section III – Computations for Scenario 2

We use the notation {G,H} to indicate that one coin comes from G and one from H; and the notations {G,G} and {H,H} to indicate both coins come from the same distribution.

The mean value m for the probability that a random coin will flip heads in Scenario 2 is given by $m = \frac{\alpha_1+\alpha_2\beta}{2(1+\beta)}$ as in Scenario 1.

## Theorem 2 (Results for Scenario 2)

$$m(U\to A) = \frac{1}{(1+\beta)^2}\left(\frac{\alpha_1}{2}\right) + \frac{\beta^2}{(1+\beta)^2}\left(\frac{\alpha_2}{2}\right) + \frac{2\beta}{(1+\beta)^2}\left(\frac{\alpha_1\alpha_2}{\alpha_1+\alpha_2}\right)$$

$$m(U\to B) = \frac{1}{(1+\beta)^2}\left(\frac{\alpha_1}{2}\right) + \frac{\beta^2}{(1+\beta)^2}\left(\frac{\alpha_2}{2}\right) + \frac{2\beta}{(1+\beta)^2}\left(\frac{\alpha_1+\alpha_2-\alpha_1\alpha_2}{4-(\alpha_1+\alpha_2)}\right)$$

**Proof:** We proceed to compute the probabilities and means associated with A and B. The mean of coins assigned from pairs of the form {G,G} is $\alpha_1/2$, the mean of G. The corresponding statement holds for coins assigned from pairs of the form {H,H}. We arrange the results as follows.

| Pair | Probability of Appearance | P(U→A) | EH(A) |
|---|---|---|---|
| {G,G} | $(1/(1+\beta))^2$ | ½ | $\alpha_1/4$ |

With probability 1, a coin from G will be examined. It will be assigned to B with probability 1/2, which is also the probability that a coin from G will be assigned to A. The expected number of heads flipped by the coin assigned to A is $\alpha_1/4$ .

| Pair | Probability of Appearance | P(U→A) | EH(A) |
|---|---|---|---|
| {H,H} | $(\beta/(1+\beta))^2$ | ½ | $\alpha_2/4$ |

With probability 1, a coin from H will be examined. It will be assigned to B with probability 1/2, which is also the probability that a coin from H will be assigned to A. The expected number of heads flipped by the coin assigned to A is $\alpha_2/4$ .

| | | | |
|---|---|---|---|
| {G,H} | $2\beta/(1+\beta)^2$ | $\frac{\alpha_1+\alpha_2}{4}$ | $\frac{\alpha_1\alpha_2}{4}$ |

With probability 1/2, a coin from G will be examined, as it is already known that the pair consists of one coin from G and one from H. It will flip heads with probability q, so the combined probability that a coin from G will be examined and U→A is q/2. The expected number of heads flipped by the coin assigned to A is pq/2.

. It will flip tails with probability 1-q, so the combined probability that a coin from G will be examined and U→B is (1-q)/2. The expected number of heads flipped by the coin assigned to B is p(1-q)/2.

With probability 1/2, a coin from H will be examined. It will flip heads with probability p, so the combined probability that a coin from G will be examined and U→A is p/2. The expected number of heads flipped by the coin assigned to A is pq/2.

. It will flip tails with probability 1-p, so the combined probability that a coin from G will be examined and U→B is (1-p)/2. The expected number of heads flipped by the coin assigned to B is q(1-p)/2.

The following three equalities are specifically for pairs of the form {G,H}.

So

$$\text{P(U}\rightarrow\text{A)} = \frac{1}{\alpha_1\alpha_2}\int_0^{\alpha_2}\int_0^{\alpha_1}\frac{1}{2}(p+q)\,dq\,dp = \frac{\alpha_1+\alpha_2}{4}$$

and

$$\text{EH(A)} = \frac{1}{\alpha_1\alpha_2}\int_0^{\alpha_2}\int_0^{\alpha_1} pq\,dq\,dp = \frac{\alpha_1\alpha_2}{4}$$

Therefore

$$\text{m(U}\rightarrow\text{A)} = \frac{EH(A)}{\text{P(U}\rightarrow\text{A)}} = \frac{\alpha_1\alpha_2}{\alpha_1+\alpha_2}$$

When all three pairs {G,G}, {H,H} and {G,H} are taken into account, we have

$$\text{m(U}\rightarrow\text{A)} = \frac{1}{(1+\beta)^2}\left(\frac{\alpha_1}{2}\right) + \frac{\beta^2}{(1+\beta)^2}\left(\frac{\alpha_2}{2}\right) + \frac{2\beta}{(1+\beta)^2}\left(\frac{\alpha_1\alpha_2}{\alpha_1+\alpha_2}\right)$$

The following three equalities are also specifically for pairs of the form {G,H}.

$$P(U\to B) = 1 - \frac{\alpha_1+\alpha_2}{4}$$

and

$$EH(B) = \frac{1}{\alpha_1\alpha_2}\int_0^{\alpha_2}\int_0^{\alpha_1}\frac{1}{2}(p(1-q)+q(1-p))\,dq\,dp = \frac{\alpha_1+\alpha_2-\alpha_1\alpha_2}{4}$$

Therefore

$$m(U\to B) = \frac{EH(B)}{P(U\to B)} = \frac{\alpha_1+\alpha_2-\alpha_1\alpha_2}{4-(\alpha_1+\alpha_2)}$$

When all three pairs {G,G}, {H,H} and {G,H} are taken into account, we have

$$m(U\to B) = \frac{1}{(1+\beta)^2}\left(\frac{\alpha_1}{2}\right) + \frac{\beta^2}{(1+\beta)^2}\left(\frac{\alpha_2}{2}\right) + \frac{2\beta}{(1+\beta)^2}\left(\frac{\alpha_1+\alpha_2-\alpha_1\alpha_2}{4-(\alpha_1+\alpha_2)}\right) \blacksquare$$

**Cor. 2.1** $m(U\to A) \le m \le m(U\to B)$, with equality occurring if and only if $\alpha_1 = \alpha_2$.

**Proof:** The proof that follows is performed using the symbol $\le$, the terminal inequality in that proof shows that equality holds if and only if $\alpha_1 = \alpha_2$.

First, put both $m(U\to A)$ and m over the common denominator $2\,(1+\beta)^2\,(\alpha_1 + \alpha_2)$. When we do so, we obtain

numerator of $m(U\to A) = (\alpha_1 + \alpha_2\beta^2)\,(\alpha_1 + \alpha_2) + 4\alpha_1\alpha_2\beta$

numerator of $m = (\alpha_1 + \alpha_2\beta)\,(\alpha_1 + \alpha_2)\,(1 + \beta)$

Showing that $m(U\to A) \le m$ now requires merely showing that the numerator of $m(U\to A) \le$ numerator of m. Expanding both and canceling the common term $\alpha_2\beta^2\,(\alpha_1 + \alpha_2)$ reduces the problem to showing $\alpha_1\,(\alpha_1 + \alpha_2) + 4\alpha_1\alpha_2\beta \le (\alpha_1 + \alpha_2\beta)\,(1 + \beta)\,(\alpha_1 + \alpha_2)$. Expanding the right side. canceling the common term $\alpha_1 (\alpha_1 + \alpha_2)$ and then canceling the common factor $\beta$ reduces the problem to showing that $4\alpha_1\alpha_2 \leq (\alpha_1 + \alpha_2)^2$. Since this is equivalent to showing $(\alpha_1 - \alpha_2)^2 \geq 0$, we have shown $m(U\to A) \leq m$.

Now, put both $m(U\to B)$ and m over the common denominator $2\,(1+\beta)^2\,(4 - (\alpha_1 + \alpha_2))$. When we do so, we obtain

numerator of $m = (\alpha_1 + \alpha_2\beta)\,(4 - (\alpha_1 + \alpha_2))\,(1 + \beta)$

numerator of $m(U\to B) = (\alpha_1 + \alpha_2\beta^2)\,(4 - (\alpha_1 + \alpha_2)) + 4\beta\,(\alpha_1 + \alpha_2 - \alpha_1\alpha_2)$

Let $c = \alpha_1 + \alpha_2$; we want to show $(\alpha_1 + \alpha_2\beta)\,(4 - c)\,(1 + \beta) \leq (\alpha_1 + \alpha_2\beta^2)\,(4 - c) + 4\beta\,(c - \alpha_1\alpha_2)$. This is equivalent to showing $0 \leq (4 - c)(\alpha_1 + \alpha_2\beta^2 - (\alpha_1 + \alpha_2\beta)\,(1 + \beta)) + 4\beta\,(c - \alpha_1\alpha_2)$. Expanding the factor of $(4 - c)$, canceling and collecting terms shows that this is equivalent to showing $0 \leq (4 - c)(-\beta c) + 4\beta\,(c - \alpha_1\alpha_2)$. This is equivalent to showing $0 < c^2 - 4\alpha_1\alpha_2$; and that was how the proof that $m(U\to A) \leq m$ concluded. ∎